\documentclass{siamart251216}
\usepackage{graphicx} 
\usepackage{tikz}
\usetikzlibrary{arrows.meta}
\usepackage{comment}
\definecolor{meshink}{RGB}{103,120,133}
\definecolor{boundaryink}{RGB}{35,51,65}
\definecolor{capfill}{RGB}{219,237,248}
\definecolor{capink}{RGB}{20,103,151}
\definecolor{obstaclefill}{RGB}{235,237,240}

\usepackage{microtype}
\usepackage{ragged2e}

\usepackage{amssymb,mathtools,dsfont}
\numberwithin{equation}{section}

\newsiamremark{remark}{Remark}

\newcommand{\tfa}{\text{ for all }}
\newcommand{\tfor}{\text{ for }}

\newcommand{\R}{\mathbb R}

\newcommand{\cT}{\mathcal T}

\renewcommand{\Im}{\operatorname{Im}}

\usepackage[show]{ed}

\usepackage{soul}

\definecolor{escol}{rgb}{0,0,0.8}
\definecolor{estcol}{rgb}{0,0.5,0}
\definecolor{esnewcol}{rgb}{0,0.5,0}

\theoremstyle{plain}
\newtheorem{main theorem}{Theorem}

\makeatletter
\newcommand{\newparallel}{\mathrel{\mathpalette\new@parallel\relax}}
\newcommand{\new@parallel}[2]{%
  \begingroup
  \sbox\z@{$#1T$}
  \resizebox{!}{\ht\z@}{\raisebox{\depth}{$\m@th#1/\mkern-5mu/$}}%
  \endgroup
}
\makeatother

\title{For any wavenumber, the 
Helmholtz FEM matrix is invertible for almost all meshes}

\headers{Helmholtz matrix invertibility for almost all meshes}
  {J. Galkowski, M. Meliani, and E. A. Spence}

\author{Jeffrey Galkowski\thanks{Department of Mathematics,
  University College London, London, UK.
  (\email{j.galkowski@ucl.ac.uk}).}
\and Mostafa Meliani\thanks{Department of Mathematical Sciences,
  University of Bath, Bath, UK.
  (\email{mm4138@bath.ac.uk}, \email{eas25@bath.ac.uk}).}
\and Euan A.~Spence\footnotemark[2]}
\date{}

\begin{document}

\maketitle
\begin{abstract}
We study invertibility of the conforming piecewise-linear Galerkin matrix for Helmholtz problems with a Dirichlet obstacle and the radiation condition approximated by \emph{either} 
impedance boundary condition \emph{or} a 
complex absorbing potential. 
We prove that for almost every simplicial mesh, and every fixed wavenumber, the Helmholtz FEM  matrix with piecewise-linear elements is invertible.
\end{abstract}

\section{Introduction}

\subsection{The main result}

Let $\Omega \subset \mathbb{R}^d$ be a bounded, Lipschitz domain with $\partial\Omega=\Gamma_D\sqcup\Gamma_R$.
Define the sesquilinear form
\begin{equation}
    \begin{aligned}
 &a_k(u,v)\\ 
 &:=
 \int_\Omega
 \left(
 (Ak^{-1}\nabla  u)\cdot\overline{k^{-1}\nabla v} -\tfrac{1}{2}i(b \cdot k^{-1}\nabla u \overline{v}-u(\overline{b\cdot k^{-1}\nabla v}))
 -
 n u\overline v -iWu\overline{v}\right) dx \\
 &\quad-ik^{-1}\int_{\Gamma_R}\beta u\overline{v} \,ds
\label{e:sesqui_generic_invertibility_thm}
\end{aligned}
\end{equation}
with $A$ real, symmetric positive definite and analytic up to $\partial \Omega$, $b$ and $n $ real valued and analytic  up to $\partial \Omega$, $W\in L^\infty(\Omega;[0,\infty))$, $\beta\in L^\infty(\Gamma_R;(0,\infty))$, 
 and $\alpha>0$ such that 
 $(A(x)\xi)\cdot \overline{\xi}\geq\alpha|\xi|^2$
 for all 
$x\in\overline\Omega,\, \xi\in\mathbb R^d.$
Note that this sesquilinear form covers two different types of absorbing boundary conditions:~an impedance condition on $\Gamma_R$, and a complex absorbing potential (where $W$ will be supported in a neighbourhood of the part of $\Gamma_D$ corresponding to the truncation boundary).

Define 
$$
\mathcal{A}:=\Big\{ \mathcal{T}\,:\, \mathcal{T} \text{
is a conforming triangulation of }\Omega\text{ with }\partial\Omega \subset \bigcup_{T\in\mathcal{T}}\partial T\Big\}.
$$
The notion of a conforming triangulation is defined in Appendix \ref{sec:App_definitions}, we sometimes informally refer to a triangulation as a mesh.
The set of \emph{nodes} of a triangulation $\mathcal{T}$ is defined by 
$$
\mathcal{N}(\mathcal{T}):=\Big\{ x\in\mathbb{R}^d\,:\, x\text{ is a vertex of an element }T\in\mathcal{T}\Big\}.
$$
The \emph{connectivity matrix of $\mathcal{T}$} is defined by $\mathcal{C}(\mathcal{T})\in\mathbb{M}(M\times (d+1))$, $M :=|\mathcal T|$, $\cT:=\{T_i\}_{i=1}^M$, $N:=|\mathcal{N}(\mathcal{T})|$, $\{x_i\}_{i=1}^N:=\mathcal{N}(\mathcal{T})$
$$
\mathcal{C}_{ij}= \ell  \quad \iff \quad x_{\ell} \textrm{ is the $j$-th vertex of } T_i
$$
(where we fix an ordering of the vertices on each $T_i$, so that the $j$-th vertex is well-defined).
Let 
$$
\mathcal{V}_{\mathcal{T}}^1:=\{ u\in H^1(\Omega)\,:\, u|_{T}\in\mathbb{P}^1(T),\, u|_{ \Gamma_D}=0\}.
$$
Informally, our main result is that for almost every $\mathcal{T}\in \mathcal{A}$ the Galerkin matrix associated to $a_k(\cdot,\cdot)$ in $\mathcal{V}_{\mathcal{T}}^1$ is invertible. To state this precisely, we introduce the following  notation.

For $N>0$ and $\mathcal{C}$ a connectivity matrix, denote 
$$
\mathcal{A}(N,\mathcal{C}):=\Bigg\{ (x_1,\dots,x_N)\in\mathbb{R}^{dN}\,:\,\begin{gathered} \mathcal{T}\in \mathcal{A},\, \mathcal{T}\text{ has connectivity }\mathcal{C},\\
|\mathcal{N}(\mathcal{T})|=N,\, \{x_j\}_{j=1}^N=\mathcal{N}(\mathcal{T})\end{gathered}\Bigg\}
$$
and denote by $\mu_{N,\mathcal{C}}$ the $r-$Hausdorff measure on $\mathcal{A}(N,\mathcal{C})$, where $r$ is the Hausdorff dimension of $\mathcal{A}(N,\mathcal{C})$.

\begin{theorem}[Invertibility of the Helmholtz FEM matrix for almost every mesh]
\label{t:generic-invertibility}
Let  $k>0$, $a_k$ as in~\eqref{e:sesqui_generic_invertibility_thm}, $N>0$, $\mathcal{C}$ be a connectivity matrix, 
 and $\mathcal{I}(N, \mathcal{C}) \subset \mathcal{A}(N,\mathcal{C})$ be the set of all $\{\tilde{x}_j\}_{j=1}^N \in \mathcal{A}(N, \mathcal{C})$ such that the mesh $\tilde{\mathcal{T}}$ with nodes $\{\tilde{x}_j\}_{j=1}^N$ and connectivity $\mathcal{C}$ satisfies
\begin{equation}
u_h\in\mathcal{V}_{\tilde{\mathcal{T}}}^1 ,\qquad a_k(u_h,v_h)=0\,\tfa v_h\in \mathcal{V}_{\tilde{\mathcal{T}}}^1 \qquad \Longrightarrow\qquad  u_h=0.
\end{equation}
Then, 
\[
\mu_{N,\mathcal{C}}(\mathcal I(N, \mathcal{C})) = \mu_{N,\mathcal{C}}(\mathcal A(N, \mathcal{C})).
\]
In particular, the Galerkin matrix is invertible $\mu_{N,\mathcal{C}}$ almost everywhere.

\end{theorem}

\subsection{The context of Theorem \ref{t:generic-invertibility}}

For the FEM applied to the Helmholtz equation with general, unstructured meshes, the Galerkin solution is known to exist and be unique when $(hk)^{2p}C_{\rm sol}$ is sufficiently small, where $p$ is the polynomial degree and $C_{\rm sol}$ is the norm of the solution operator (such that $C_{\rm sol}\sim k$ for nontrapping problems).
Furthermore, under this threshold, for $k$-oscillating data, the relative error of the Galerkin solution is $k$-uniformly bounded; 
this result was proved for general Helmholtz problems in any dimension in \cite{GS3}, following the earlier work \cite{Wu2014, DuWu2015, LiWu2019, LSW2}, with related work in \cite{melenk2010convergence, MelenkSauter2011, ChNi:20, GrSa:20}
proving $k$-uniformly quasi-optimality when $(hk)^p C_{\rm sol}$ is sufficiently small.

The threshold ``$(hk)^{2p}C_{\rm sol}$ sufficiently small" was famously identified by the 1-d investigations of Ihlenburg and Babu\v{s}ka \cite{IhlenburgBabuska1995, IhlenburgBabuska1997}, with
\cite[Theorem 4]{IhlenburgBabuska1995} and \cite[Theorem 3.3]{IhlenburgBabuska1997} also proving existence and uniqueness of Galerkin solution when $hk$ is sufficiently small in the one-dimensional case.

The solvability of the Helmholtz FEM matrix was recently investigated in detail in \cite{bernkopf2023solvability} (with the follow-up paper \cite{ChSa:26} also considering the time-harmonic Maxwell equations). 
The paper 
\cite{bernkopf2023solvability} considers \eqref{e:sesqui_generic_invertibility_thm} with $A=I, b=0, n=1,$ and $W=0$ and proves results that can be grouped into the following three classes:
\begin{enumerate}
    \item The
Helmholtz FEM matrix is always invertible in 1-d \cite[Theorem 3.1]{bernkopf2023solvability} and on 2-d 
structured quadrilateral grids for $p=1,2,3$
\cite[Theorem 3.6]{bernkopf2023solvability}.
\item An algorithm for determining when the Helmholtz FEM matrix
with piecewise-linear elements on triangular meshes is invertible \cite[Algorithm 1]{bernkopf2023solvability}. This ``marching of the zeros" algorithm (which can be understood as a discrete analogue of the unique continuation principle) propagates zeros from the boundary under geometric conditions on the mesh
(that, for certain interior edges, the sum of the two angles opposite that edge in the adjacent triangles does not exceed $\pi$).
\item There is a 2-d triangulation for which Helmholtz FEM matrix is \emph{not} invertible when $kh=c$ for some specific constant $c$ \cite[Lemma 3.3]{bernkopf2023solvability}.
\end{enumerate}
The last result shows that, for a given $k>0$, one cannot expect the Helmholtz FEM matrix to be invertible for \emph{all} meshes.
Theorem \ref{t:generic-invertibility} shows that the FEM matrix is nevertheless invertible for \emph{almost every} mesh (in the sense made precise in the statement of the result).

\subsection{The idea of the proof, and structure of the paper}\label{s:idea}

The key observation used in the proof of Theorem~\ref{t:generic-invertibility} is that, in the presence of analytic coefficients, the Galerkin matrix is an analytic function of node positions (see \S\ref{s:analyticity}). Thus, if all coefficients were analytic, it would be sufficient to exhibit a single mesh (in every connected component of admissible meshes) whose Galerkin matrix is invertible. 

The presence of a CAP in~\eqref{e:sesqui_generic_invertibility_thm} means that not all coefficients can be analytic. However, because $W\geq 0$ and $\beta> 0$, any solution, $u_h\in\mathcal{V}_{\mathcal{T}}^1$ to
$$
a_k(u_h,v_h)=0,\qquad \text{for all }v_h\in\mathcal{V}_{\mathcal{T}}^1
$$
satisfies $Wu_h=0$ and $\beta u_h|_{\Gamma_R}=0$ and hence $u_h\in H_0^1(\Omega)$ and 
$$
(\Re a_k)(u_h,v_h)=0,\qquad  \text{for all }v_h\in\mathcal{V}_{\mathcal{T}}^1,
$$
where 
$$
(\Re a_k)(u,v) := \frac{a_k(u,v) + \overline{a_k(v,u)}}{2}.
$$
In particular, $\Re a_k$ is of the form~\eqref{e:sesqui_generic_invertibility_thm} with $W=0$ and $\Gamma_R=\emptyset$. Thus, if the Galerkin matrix for $\Re a_k$ (on $H_0^1(\Omega)\cap \mathcal{V}_{\mathcal T}^1$) is invertible, then so is the Galerkin matrix for $a_k$ on $\mathcal{V}_{\mathcal{T}}^1$. 

This deduction shows that it is enough to consider a sesquilinear form with analytic coefficients and Dirichlet boundary conditions and hence one that is  analytic in all interior node positions.

\begin{figure}
    \begin{center}
%
%
\begingroup
\definecolor{amInk}{HTML}{263643}
\definecolor{amTeal}{HTML}{007D84}
\definecolor{amInactive}{HTML}{EDF2F5}
\definecolor{amActive}{HTML}{F6CFAC}
\definecolor{amOrange}{HTML}{B75113}
\begin{tikzpicture}[
  x=1.8cm,y=1.8cm,
  font=\sffamily,
  am edge/.style={draw=amInk,line width=0.38pt,line join=round},
  am boundary/.style={draw=amInk,line width=0.72pt,line join=round},
  am direction/.style={
    draw=amTeal,line width=0.65pt,
    -{Stealth[length=2.6pt,width=2.6pt]}
  }
]
  \foreach \meshTime [count=\panel from 0] in {0,0.5,0.9} {
    \begin{scope}[
      shift={({2.441176470588*\panel},0)},
      name prefix=am\panel-
    ]
      \coordinate (v0) at ( 1, 1);
      \coordinate (v1) at (-1, 1);
      \coordinate (v2) at (-1,-1);
      \coordinate (v3) at ( 1,-1);

      \foreach \j in {0,...,7} {
        \pgfmathtruncatemacro{\idx}{4+\j}
        \coordinate (v\idx) at
          ({(1/2)*cos(45*\j)},{(1/2)*sin(45*\j)});
      }

      \foreach \j/\weight in {
        0/0.2,
        1/(1-0.4*sqrt(2)),
        2/1,
        3/0.3,
        4/0.3,
        5/0.3,
        6/(0.3*(1+sqrt(2))),
        7/(1-0.4*sqrt(2))
      } {
        \pgfmathtruncatemacro{\idx}{12+\j}
        \coordinate (v\idx) at
          ({(4/5)*cos(45*\j)+\meshTime*(\weight)},
           {(4/5)*sin(45*\j)});
      }

      \foreach \a/\b/\c/\meshFill in {
      12/0/13/amActive,
      13/0/14/amActive,
      14/0/1/amInactive,
      14/1/15/amInactive,
      15/1/16/amInactive,
      16/1/2/amInactive,
      16/2/17/amInactive,
      17/2/18/amInactive,
      18/2/3/amInactive,
      18/3/19/amInactive,
      19/3/12/amActive,
      12/3/0/amActive,
      4/12/13/amInactive,
      4/13/5/amInactive,
      5/13/14/amInactive,
      5/14/6/amInactive,
      6/14/15/amInactive,
      6/15/7/amInactive,
      7/15/16/amActive,
      7/16/8/amActive,
      8/16/17/amActive,
      8/17/9/amActive,
      9/17/18/amInactive,
      9/18/10/amActive,
      10/18/19/amInactive,
      10/19/11/amInactive,
      11/19/12/amInactive,
      11/12/4/amInactive
      } {
        \fill[\meshFill] (v\a)--(v\b)--(v\c)--cycle;
      }

      \foreach \a/\b in {
      0/1,0/3,0/12,0/13,0/14,1/2,1/14,1/15,%
      1/16,2/3,2/16,2/17,2/18,3/12,3/18,3/19,%
      4/5,4/11,4/12,4/13,5/6,5/13,5/14,6/7,%
      6/14,6/15,7/8,7/15,7/16,8/9,8/16,8/17,%
      9/10,9/17,9/18,10/11,10/18,10/19,11/12,11/19,%
      12/13,12/19,13/14,14/15,15/16,16/17,17/18,18/19%
      } {
        \draw[am edge] (v\a)--(v\b);
      }
      \draw[am boundary] (v0)--(v1)--(v2)--(v3)--cycle;
      \draw[am boundary]
        (v4)--(v5)--(v6)--(v7)--(v8)--(v9)--(v10)--(v11)--cycle;

      \foreach \j in {0,...,11}
        \fill[amInk] (v\j) circle[radius=0.97pt];
      \foreach \j in {12,...,19}
        \fill[amTeal] (v\j) circle[radius=1.5pt];

      \ifnum\panel=0\relax
        \foreach \j in {0,...,7} {
          \pgfmathtruncatemacro{\idx}{12+\j}
          \path (v\idx) ++({45*\j}:0.12)
            node[text=amTeal,font=\scriptsize,inner sep=0pt] {};
        }
      \fi

      \node[anchor=base,text=amInk,font=\bfseries\boldmath]
        at (0,1.147059) {$t=\meshTime$};
      \draw[am direction] (-0.205882,-1.235294)--(0.117647,-1.235294);
      \node[anchor=base west,text=amTeal,font=\small,inner sep=0pt]
        at (0.176471,-1.270588) {$e$};
    \end{scope}
  }

  \filldraw[fill=amActive,draw=amOrange,line width=0.3pt]
    (1.06,-1.49) rectangle (1.17,-1.40);
  \node[anchor=west,text=amInk,font=\small,inner sep=0pt]
    at (1.24,-1.445) {Collapsing triangles};
  \fill[amTeal] (3.15,-1.445) circle[radius=1.5pt];
  \node[anchor=west,text=amInk,font=\small,inner sep=0pt]
    at (3.24,-1.445) {Moving vertices};
\end{tikzpicture}
\endgroup
    \end{center}
    \caption{An example of the mesh deformation used to prove Theorem~\ref{t:generic-invertibility}
    }
    \label{fig:collapse}
\end{figure}
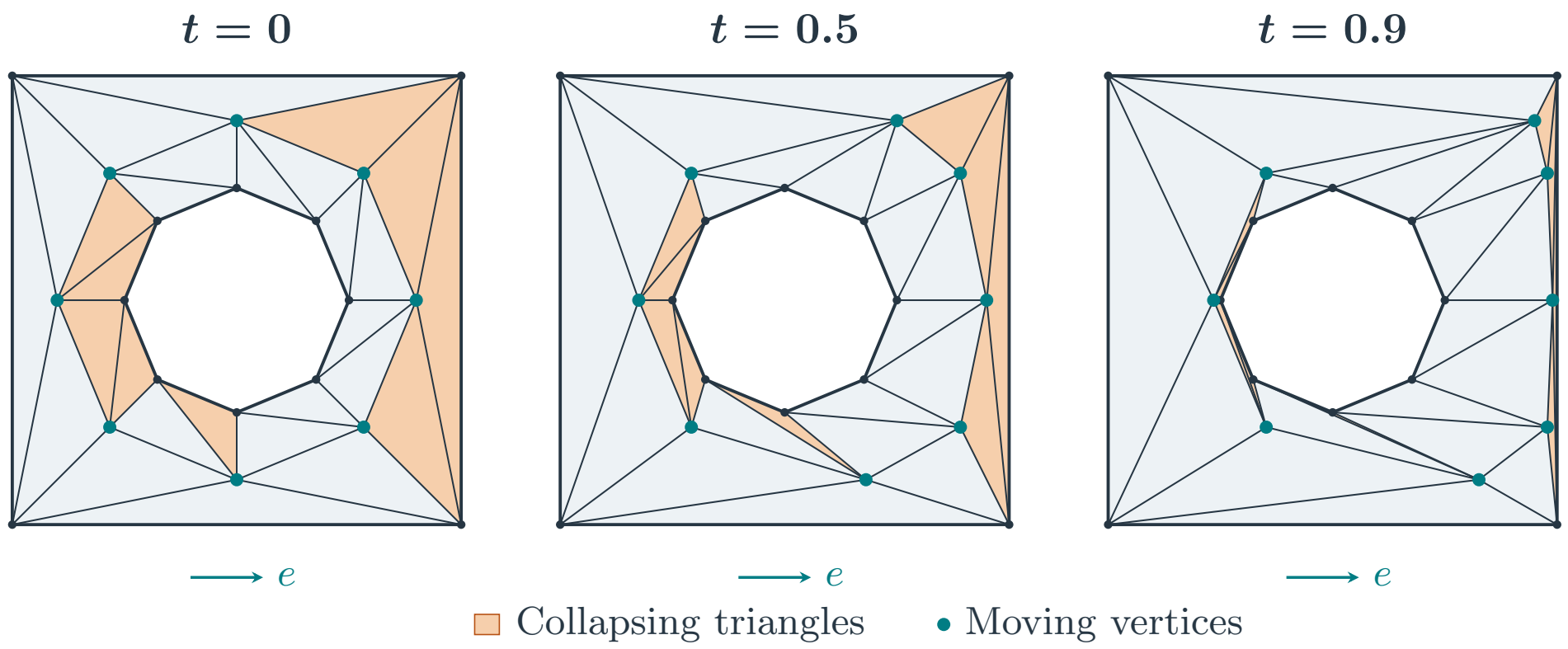

We then proceed to exhibit a mesh for which the Galerkin matrix for $\Re a_k$ is invertible, via proving that $\Re a_k$ is coercive. To obtain coercivity, we use three ideas
\begin{enumerate}
    \item If we deform the mesh so that mesh elements collapse, the derivatives of piecewise-linear basis functions are likely to grow compared to their nodal values. 
    Since the boundary is fixed, it is not possible to collapse all mesh elements, but we would like to collapse a maximal (in an appropriate sense) set of mesh elements. 
    \item One way to collapse elements is to push nodes in a single direction and, since boundary nodes need to be fixed, this deformation takes the form
   \begin{equation}
       \label{e:prelimDeform}
       x\mapsto x+tw(x)e
   \end{equation} for some unit vector $e$ and $w\in\mathcal{V}_{\mathcal{T}}^1\cap H_0^1(\Omega)$ with real values. This is a well defined homeomorphism  on $\Omega$ for $t\in[0,1)$ if $\partial_e w\geq -1$ and an element then collapses at $t=1$ if $\partial_ew =-1$. See Figure~\ref{fig:collapse} for an example of this type of deformation.
   \item
    We are therefore led to consider the closed, bounded, and convex set of deformations satisfying the linear constraint $\partial_e w\geq -1$ (so that the map \eqref{e:prelimDeform} remains a homeomorphism for $t\in [0,1)$).
    The intuition from linear programming (see, e.g., \cite[Chapter 2]{bertsimas1997introduction}) is that, given an extreme point $w_*$, any  $v\in\mathcal{V}_{\mathcal{T}}^1\cap H_0^1(\Omega)$ is determined by the values of the constraint function, $\partial_e v$, on the elements `activated' by $w_*$ (this is made precise in Lemma \ref{l:PFType} below); here 
the constraint on an element $T$ is `active' if and only if the extreme point, $w_*$, satisfies $\partial_e w_*=-1$ and hence if and only if the element collapses at $t=1$. 
The active elements in Figure \ref{fig:collapse} are shown in orange.
\end{enumerate}
The precise description of the mesh deformation obtained following these  ideas is contained in \S\ref{s:deform} and 
\S\ref{s:choice}.

Finally, \S\ref{s:proof} contains the proof of the main theorem. The fact that the active constraints control the function leads to a Poincar\'e-type inequality (see Lemma \ref{l:PFType}) that controls a finite element function by its derivative on a subset, $\mathcal{E}$, of elements that collapse as $t\to 1^-$. This collapse in turn then shows that the constant in the Poincar\'e inequality becomes arbitrarily small as $t\to 1^-$ and hence that $\Re a_k$ is coercive on $\mathcal{V}_{\mathcal{T}_t}^1\cap H_0^1(\Omega)$ for $t$ close to 1.

\subsection{Use of AI}

During the development of this article, the authors identified the
following strategy:~for each fixed wavenumber, construct an admissible mesh on which the discrete Helmholtz problem is invertible, and then use the analytic dependence of the assembled matrix on the vertex coordinates to establish invertibility for almost every admissible mesh. The main difficulty was to construct a deformation that preserves conformity, connectivity, and element orientation while producing the required coercivity.

The authors then used ChatGPT 6 Astra  as a brainstorming aid in developing the mesh deformation. After the authors suggested that it might be possible to make many elements collapse so that the effective Poincar\'e constant on finite elements becomes small, ChatGPT proposed the idea of using the form of mesh deformation~\eqref{e:prelimDeform} with an extremal $w$. The authors then realized that this comes from a general principle in linear programming.  ChatGPT was also used to generate Figure~\ref{fig:collapse}.

All text in the article was written by the authors, and the authors assume responsibility for all content.

\section{Analyticity of the Galerkin matrix}\label{s:analyticity}

In this section, we show that the Galerkin matrix for $\Re a_k(\cdot, \cdot)$ is analytic in certain nodes of the mesh. For this, we define the \emph{Galerkin matrices of $a_k(\cdot, \cdot)$ and $\Re a_k(\cdot,\cdot)$ with triangulation $\mathcal{T}$} by 
$$
\big(G(a_k,\mathcal{T})\big)_{ij}:=  a_k(\tilde{\phi}_i,\tilde{\phi}_j),\qquad \big(G(\Re a_k,\mathcal{T})\big)_{ij}:=(\Re a_k)(\phi_i,\phi_j)
$$
where $\{\tilde{\phi}_{i}\}_{i=1}^{\tilde{M}}$ is the nodal basis of $\mathcal{V}_{\mathcal{T}}^1$, $\{\phi_i\}_{i=1}^{M}$ is the nodal basis for $\mathcal{V}_{\mathcal{T}}^1\cap H^1_0(\Omega)$.

We also define
\begin{gather*}
\mathcal{M}_{\mathcal{T}}:=\mathcal{N}(\mathcal{T})\setminus\partial\Omega,\qquad \mathcal{Z}_{\mathcal{T}}:=\mathcal{N}(\mathcal{T})\setminus \mathcal{M}_{\mathcal{T}}.
\end{gather*}

\begin{lemma}
\label{l:analytic}
    Let $\mathcal{T}\in \mathcal{A}$ with connectivity $\mathcal{C}$ and $\{x_j\}_{j=1}^J=\mathcal{M}_{\mathcal{T}}$. Then, there is a neighborhood $U$ of $(x_1,\dots,x_J)$ such that, with $y\in\big(\mathbb{R}^{d}\big)^J$ and
    $
    \mathcal{T}_y
    $ the triangulation of $\Omega$ defined by moving $x_j$ to $y_j$ for $j=1,\dots, J$ that has connectivity $\mathcal{C}$, the map
    $$
    U\ni y\mapsto G(\Re a_k,\mathcal{T}_y)\in\mathbb{M}(M\times M)
    $$
    is analytic.
\end{lemma}

\begin{proof}[Proof of Lemma \ref{l:analytic}]
Let 
$$
\mathcal{C}_{j}:=\{ x\in\mathcal{N}(\mathcal{T})\,:\, \exists T\in\mathcal{T}\text{ with }x\in T\text{ and }x_j\in T\}.
$$
It is enough to show that for any $j=1,\dots, J$ 
$$
\big(G(\Re a_k,\mathcal{T}_y)\big)_{il}
,\qquad x_i,x_l\in\mathcal{C}_j
$$
is analytic in $y$ since all other entries are independent of $y$. 

Fix the reference element 
$$
\hat{T}:=\Big\{ 0\leq x_l,\, \sum_{l=1}^dx_l\leq 1\Big\}
$$
and the reference nodal basis $\{\widehat\phi_l\}_{l=0}^d$,
\[
 \widehat\phi_0(x)
 =
 1-\sum_{r=1}^d\widehat x_r,
 \qquad
 \widehat\phi_l(\widehat x)= x_l
 \quad(1\leq l\leq d).
\]
Then, for a simplex $T:=\text{convex hull}(\{ q_l\}_{l=0}^{d})$ define the affine map
\begin{equation}
 F_T(x)
 :=
 q_0+B_T x,
 \qquad
 B_T
 :=
 \begin{bmatrix}
 q_1-q_0&\cdots&q_d-q_0
 \end{bmatrix}\in \mathbb{M}(d\times d).
 \label{e:element-affine-map}
\end{equation}
and the nodal basis functions $\phi_{l,T}\circ F_T:=\hat{\phi}_l$. 

We compute the stiffness matrix. The other terms are handled similarly. 
 Since $\hat{\phi}_l$ are real, the local stiffness entries satisfy
\begin{align}
 S_{l_1l_2}(y)
 &:=
 \int_{T_y}
  A\nabla\phi_{l_1,T_y}(x)\cdot \nabla\phi_{l_2,T_y}(x) dx
 \notag\\
 &=
 |\det B_{T_y}|
 \int_{\widehat T} \Big( A(F_{T_y}(x))
 B_{T_y}^{-\top}\nabla  \hat{\phi}_{l_1}(x)\Big)
 \cdot
 \bigl(B_{T_y}^{-\top}\nabla \hat{\phi}_{l_2}(x)\bigr)
 \,d x
 \label{e:rational-local-stiffness}
\end{align}
Since $B_{T_y}$ has real coefficients, $y\mapsto B_{T_y}$ is analytic, and $B_{T_{\vec{x}}
}$ is invertible (with $\vec{x}=(x_1,\dots, x_J)$), $B_{T_y}^{-1}$, $|\det B_{T_y}|$, and $F_{T_y}(x)$ are analytic in $y$. 
Since $A$ is analytic, this completes the proof. 
\end{proof}

\section{Properties of a $1$-dimensional mesh deformation}
\label{s:deform}
Fix $e\in\mathbb{S}^{d-1}$ and $w\in \mathcal{K}_{\mathcal{T}, \R}^1$, where
\begin{equation}
 \mathcal{K}_{\mathcal{T}, \R}^1:=\big\{u\in \mathcal{V}_{\mathcal{T}}^1\,:\, u(x_j)=0 \tfor x_j\in\mathcal{Z}_{\mathcal{T}},\, u \text{ is real valued}\big\},
 \label{e:mesh-space}
\end{equation}
and define
\begin{equation} \label{e:deformation}
F(t,w,e)(x):= x+ tw(x)e
\end{equation}

\begin{remark}
    Note that $\mathcal{K}_{\mathcal{T}, \R}^1$ consists of those functions in $\mathcal{V}_{\mathcal{T}}^1\cap H_0^1(\Omega)$ that are real valued.
\end{remark}

We need a few elementary properties of $F$. 
\begin{lemma}
\label{l:goodStuff}
For all $w\in\mathcal{K}_{\mathcal{T},\R}^1$ with 
$$
(1+\partial_e w)|_{T}\geq 0,\qquad T\in\mathcal{T},
$$
and all $t\in [0,1)$ the following are true.
    \begin{enumerate}[(i)]
        \item $F(t,w,e)(\mathcal{Z}_{\mathcal{T}})=\mathcal{Z}_{\mathcal{T}}$
        \item For all $T\in\mathcal{T}$, 
        \begin{equation}\label{e:determinant}
        \det D F(t,w,e)(x)= (1+t\partial_e w(x)),\qquad x\in T.
        \end{equation}
        \item For all $T\in \mathcal{T}$ and $u\in H^1(F(t,w,e)(T))$
        \begin{equation}
        \label{e:derivativeChange}
         \|\partial_e (u\circ F(t,w,e))\|^2_{L^2(T)}= \big|(1+t\partial _e w)|_{T}\big|\|\partial_e u\|^2_{L^2(F(t,w,e)(T))}
        \end{equation}
    \end{enumerate}
\end{lemma}
\begin{proof}
    (i) This follows from the definition of $\mathcal{K}_{\mathcal{T},\mathbb{R}}^1$. 

    (ii) By direct calculation,
    $$
    DF= I+t e\otimes \nabla w(x).
    $$
    Hence, since $\det(1+ u\otimes v)= 1 + u\cdot v$,
    $$
    \det DF= (1+t\partial_e w).
    $$

    (iii) By the chain rule $\nabla (u\circ F)(x) = (DF(x))^T \nabla u(F(x))$, so that
    $$
    \partial _e (u\circ F)= \langle (\nabla u)\circ F,\partial_e F\rangle = \langle (\nabla u)\circ F, (1+t\partial_e w)e\rangle = (1+t\partial_e w)(\partial_e u)\circ F.
    $$
    The equality then follows by changing variables, using the fact that $\det DF\neq 0$. 
\end{proof}

\begin{lemma}
\label{l:reallyGoodStuff}
For all $w\in\mathcal{K}_{\mathcal{T},\R}^1$ with
\begin{equation}\label{e:chicken1}
(1+\partial_e w)|_{T}\geq 0,\qquad T\in\mathcal{T},
\end{equation}
and all $t\in [0,1)$, the following are true.
    \begin{enumerate}[(i)]
     \item $F(t,w,e):\overline{\Omega}\to \overline{\Omega}$ is a homeomorphism.
        \item $\mathcal{T}_t:=\{ F(t,w,e)(T)\}_{T\in\mathcal{T}}$ is a triangulation of $\Omega$ in $\mathcal{A}$ with the same connectivity as $\mathcal{T}$.
        \item $F(t,w,e)(\mathcal{M}_{\mathcal{T}})=\mathcal{M}_{\mathcal{T}_t}$.
        \end{enumerate}
\end{lemma}
\begin{proof}
(i) Without loss of generality, we assume $e=(1,0,\dots,0)$. Then,
$$
F(t,w,e)(x_1,x')=(x_1+tw(x_1,x'),x').
$$
Since $F$ preserves $x'$ it is therefore enough to show that for each fixed $x_0'\in\mathbb{R}^{d-1}$ 
$$
\tilde{F}(x_1):= F(t,w,e)(x_1,x_0')
$$
is a homeomorphism from $\Omega_{x_0'}:= \Omega\cap\{x'=x'_0\}$ to itself.

Since $0\leq t<1$ and $1+\partial_e w\geq 0$, $\tilde{F}$ is strictly increasing in $x_1$ and hence is a homeomorphism onto its image. To see that $\tilde{F}:\Omega_{x_0'}\to \Omega_{x_0'}$ is surjective, observe that if $y\in \partial \Omega_{x_0'}$, then $w(y,x_0')=0$ and hence $\tilde{F}(y)=y$. Since $\Omega_{x_0'}$ is open it is a countable union of disjoint open intervals and hence, by the intermediate value theorem on each interval, $\tilde{F}$ is surjective and hence a homeomorphism.

 (ii) By~\eqref{e:determinant} and \eqref{e:chicken1}, on each element $T\in\mathcal{T}$ $F(t,w,e)$ is a non-degenerate, orientation-preserving affine map and hence the fact that $\mathcal{T}_t$ is a triangulation follows immediately; that its connectivity is the same as that $\mathcal{T}$ follows from the fact that $F(t,w,e)$ is a homeomorphism. 

(iii) This follows from the definition of $\mathcal{K}_{\mathcal{T},\mathbb{R}}^1$.    
\end{proof}

\section{Choice of the mesh deformation and linear programming}\label{s:choice}

We now choose $w$ in the deformation~\eqref{e:deformation} to collapse enough elements so that a finite element function is controlled by its derivative. 
To do this, we consider the following set
$$
\mathcal{P}_e:=\big\{ w\in\mathcal{K}_{\mathcal{T},\mathbb{R}}^1\,:\, (1+\partial_e w)|_{T}\geq 0,\, T\in\mathcal{T}\big\}
$$ 
To maximize the number of collapsing elements of the mesh (i.e., elements for which $1+\partial_ew(x)=0$), we choose an extreme point of $\mathcal{P}_e$. 
(Recall $x\in K$ is an extreme point if $x=\lambda y + (1-\lambda)z$ for $y,z\in K$ and $\lambda\in (0,1)$ implies $x=y=z$.)

\begin{lemma}
    The set $\mathcal{P}_e$ is a closed, convex, bounded set. In particular, it has extreme points.
\end{lemma}
\begin{proof}
    Convexity and closedness of $\mathcal{P}_e$ follow directly from the definition. 

   Since $\mathcal{K}_{\mathcal{T},\mathbb{R}}^1$ is finite dimensional, we can show boundedness in any norm. We claim that 
   $$
   \mathcal{P}_e\subset \{ u\in C^0(\overline\Omega)\,:\, \|u\|_{L^\infty}\leq \text{diam}(\Omega)\}.
   $$
   Suppose by contradiction that $w\in\mathcal{P}_e$ with $\|w\|_{\infty}>\text{diam}(\Omega) $. Then, since $w\in C^0$, there is $x_0\in\Omega$ with $|w(x_0)|> \text{diam}(\Omega)$. We assume $w(x_0)> \text{diam}(\Omega)$; the other case is similar. Now define
    $$
    y:=\inf\{ z>0\,:\, x_0+e z\notin \Omega\}<\text{diam}(\Omega)
    $$
    and observe that $w(x_0+ye)=0$.
The derivative condition in $\mathcal{P}_e$ implies that 
$w(x_0)\leq w(x_0+y e)+y$ so that 
    $$
    \text{diam}(\Omega)<w(x_0)\leq w(x_0+y e)+y\leq \text{diam}\Omega,
    $$
    which is a contradiction. 
\end{proof}

The following Poincar\'e-type inequality is the analogue in our situation of the standard result from linear programming that, at an extreme point of a linear constraint on a bounded, convex set, the active constraints  control the vector. 
\begin{lemma}
\label{l:PFType}
    Suppose that $w_*$ is an extreme point of $\mathcal{P}_e$. Then there is $C>0$ such that for all $v\in\mathcal{K}_{\mathcal{T}}^1:=\mathcal{K}_{\mathcal{T},\mathbb{R}}^1+i\mathcal{K}_{\mathcal{T},\mathbb{R}}^1= \mathcal{V}_{\mathcal{T}}^1\cap H_0^1(\Omega)$, 
    $$
    \|v\|_{L^2}\leq C\sum_{T\in\mathcal{E}}\|\partial_e v\|_{L^2(T)},\qquad \mathcal{E}:=\big\{ T\in\mathcal{T}\,:\, 1+\partial_e w_*=0\big\}.
    $$
\end{lemma}
\begin{proof}
    Since $\mathcal{K}_{\mathcal{T}}^1$ is finite dimensional, it is sufficient to show that the map $P:\mathcal{K}_{\mathcal{T}}^1 \to \oplus_{T\in\mathcal{E}}L^2(T)$, 
    $$
    P v:= (\partial_e v|_{T})_{T\in\mathcal{E}}
    $$
    is injective. 
    Suppose that $v\in \mathcal{K}_{\mathcal{T}}^1$ is such that $\partial_e v|_{T}=0$ for all $T\in\mathcal{E}$. Then, $\Re v\in\mathcal{K}_{\mathcal{T},\mathbb{R}}^1$ and for all $s\in\mathbb{R}$
    and $T\in \mathcal{E}$
    $$
    (1+\partial_e (w_*+s\Re v))|_{T}=0.
    $$
    Moreover, since $\mathcal{T}$ is finite and  $1+\partial_e w_*>0$ for $T\notin \mathcal{E}$, there is $\delta>0$ such that 
    $$
    1+\partial_e (w_*+s\Re v)|_{T}>0,\qquad T\notin \mathcal{E},\, |s|<\delta.
    $$
    Thus $w_*+s\Re v\in\mathcal{P}_e$.
 However, since $w_*$ is an extreme point of the convex set $\mathcal{P}_e$, this implies that $\Re v=0$. An identical argument shows that $\Im v=0$; thus $v=0$ and $P$ is injective as claimed.
\end{proof}

\section{Proof of the main theorem}
\label{s:proof}

The next lemma combines Lemma~\ref{l:PFType} with the fact that elements in $\mathcal{E}$ are collapsing as $t\to 1^-$ to show that, for $u\in \mathcal{K}_{\mathcal{T}_t}^1$, an arbitrarily-small multiple of the gradient of $u$ controls the $L^2$ norm of $u$.
\begin{lemma}
\label{l:collapse}
    Let $e\in\mathbb{S}^{d-1}$ and suppose that $w_*$ is an extreme point of $\mathcal{P}_e$. Then there is $C>0$ such that, defining $\mathcal{T}_t:=F(t,w_*,e)(\mathcal{T})$,  for  all $u\in\mathcal{K}_{\mathcal{T}_t}^1$, 
    $$
    \|u\|_{L^2}^2\leq C(1-t)\|\nabla u\|_{L^2}^2.
    $$
\end{lemma}
\begin{proof}
First observe that, with $F_t:=F(t,w_*,e)$, $u\circ F_t\in\mathcal{K}_{\mathcal{T}}^1$. Hence, by~\eqref{e:derivativeChange}, 
    \begin{align*}
    \|\nabla u\|_{L^2}^2=\sum_{T\in\mathcal{T}}\|\nabla u\|_{L^2(F_t(T))}^2
    &\geq \sum_{T\in\mathcal{T}}\|\partial_e u\|_{L^2(F_t(T))}^2\\
    &=  \sum_{T\in\mathcal{T}}\big|(1+t\partial_e w_*)|_{T}\big|^{-1}\|\partial_e (u\circ F_t)\|_{L^2(T)}^2\\
    &\geq  (1-t)^{-1}\sum_{T\in\mathcal{E}}\|\partial_e (u\circ F_t)\|_{L^2(T)}^2
    \end{align*}
    Then, by Lemma~\ref{l:PFType}  and~\eqref{e:determinant} (to change variables), 
    \begin{align*}
    \|\nabla u\|_{L^2}^2&\geq   c(1-t)^{-1}\sum_{T\in\mathcal{T}}\| (u\circ F_t)\|_{L^2(T)}^2\\
&\geq c(1-t)^{-1}\|u \|^2_{L^2},
 \end{align*}
 as claimed.
\end{proof}

\begin{proof}[Proof of Theorem~\ref{t:generic-invertibility}]
Let $\mathcal{T}_0\in\mathcal{A}$ and $\mathcal{C}$ be its connectivity matrix and
    $$
    \mathcal{A}_{\mathcal{C},\mathcal{T}_0}:=\{ \mathcal{T}\in\mathcal{A}\,:\, \mathcal{T} \text{ has connectivity }\mathcal{C},\, \mathcal{Z}_{\mathcal{T}}=\mathcal{Z}_{\mathcal{T}_0}\}.
    $$
    We show that there is $W_{\mathcal{C},\mathcal{T}_0}\subset \mathcal{A}_{\mathcal{C},\mathcal{T}_0}$ open, dense, and with full measure such that 
    \begin{equation}
    \label{e:claimStep1}
        \text{ $\mathcal{T}\in W_{\mathcal{C},\mathcal{T}_0}$ implies $G(a_k,\mathcal{T})$ is invertible.}
    \end{equation}

First, suppose that 
$
G(a_k,\mathcal{T})
$
is not invertible.
Then there is $0\neq v_h\in\mathcal{V}_{\mathcal{T}}^1$ such that 
$$
a_k(v_h,w_h)=0,\qquad w_h\in\mathcal{V}_{\mathcal{T}}^1.
$$
In particular
    $$
    0=\Im a_k(v_h,v_h)= -\int_{\Omega} W|v_h|^2dx -k^{-1}\int_{\Gamma_R}\beta |v_h|^2ds,
    $$
    so that $Wv_h=0$ almost everywhere and $v_h\in H_0^1(\Omega)$. This implies that 
    $$
    (\Re a_k)(v_h,w_h)=0,\qquad \tfa  w_h\in\mathcal{V}_{\mathcal{T}}^1.
    $$
Thus, $G(a_k,\mathcal{T})$ is not invertible implies $G(\Re a_k,\mathcal{T})$ is not invertible. Now
$$
\{\mathcal{T}: G(\Re a_k,\mathcal{T})\text{ is not invertible}\}=\{\mathcal{T}: \det G(\Re a_k,\mathcal{T})=0\},
    $$
    and, by Lemma~\ref{l:analytic}, $G(\Re (a_k),\mathcal{T})$ is analytic in $\mathcal{A}_{\mathcal{C},\mathcal{T}_0}$. Therefore~\cite[Proposition 0]{Mityagin2020} shows that 
    it is sufficient to find $\mathcal{T}'\in \mathcal{A}_{\mathcal{C},\mathcal{T}_0}$ such that $G(\Re a_k,\mathcal{T}')$ is invertible. 

    Let $e\in\mathbb{S}^{d-1}$ and $w_*$ be an extreme point of $\mathcal{P}_e(\mathcal{T}_0)$. Then,  by Lemmas~\ref{l:goodStuff} and~\ref{l:reallyGoodStuff}, $\mathcal{T}_t:= F(t,w_*,e)(\mathcal{T}_0)\in \mathcal{A}_{\mathcal{C},\mathcal{T}_0}$. Now, suppose that $v_h\in\mathcal{V}_{\mathcal{T}_t}^1\cap H_0^1(\Omega)$ is such that 
    $$
    (\Re a_k)(v_h,v_h)=0.
    $$
 Then, by Lemma~\ref{l:collapse} and the G\aa rding inequality for $a_k$,
    \begin{align*}
    0=\Re a_k(v_h,v_h)\geq ck^{-2}\|\nabla v_h\|_{L^2}^2-C\|v_h\|_{L^2}^2
    \geq ck^{-2}(1-t)^{-1}\| v_h\|_{L^2}^2-C\|v_h\|_{L^2}^2.
    \end{align*}
    If $(1-t)<c(2Ck^2)^{-1}$ then  $\|v_h\|_{L^2}=0$ and hence $v_h=0$. 
In other words, for $t$ close enough to $1$, $G(\Re(a_k),\mathcal{T}_t)$ is invertible; this  proves~\eqref{e:claimStep1}, since $\mathcal{T}_0$ with connectivity $\mathcal{C}$ is arbitrary.

    To finish the proof, we now need to show that
     $$
     \mathcal{B}_{\mathcal{C}}:=\{ \mathcal{T}\in\mathcal{A}\,:\, \mathcal{T}\text{ has connectivity } \mathcal{C},\, G(a_k,\mathcal{T})\text{ is not invertible}\}
    $$
     has zero measure. 
     To do this, 
     we first express $\mathcal{B}_{\mathcal{C}}$ as a finite union of sets of meshes with a fixed number of  boundary nodes.  I.e., 
     let $L\in \mathbb{Z}^+$ and define
     \begin{gather*}
\mathcal{B}_{\mathcal{C},L}:=\mathcal{B}_{\mathcal{C}}\cap \mathcal{D}_{\mathcal{C},L},\\
 \mathcal{D}_{\mathcal{C},L}:=\big\{ \mathcal{T}\,:\, |\mathcal{Z}_{\mathcal{T}}|=L,\, \mathcal{T}\text{ has connectivity }\mathcal{C}\big\}.
     \end{gather*}
 It is therefore enough to show that $\mathcal{B}_{\mathcal{C},L}$ has measure zero for each $L$.

 For any set of boundary nodes, the interior vertices of the mesh can be varied in an open set. Furthermore, the fact that $\mathcal{T}$ is admissible requires the triangulation to resolve $\partial \Omega$, which fixes certain boundary nodes to higher codimension submanifolds-- resulting in a less than $|\mathcal{N}(\mathcal{T})|^d$ dimensional subset of $\mathbb{R}^{d|\mathcal{N}(\mathcal{T})|}$ However, for any admissible set of nodes, the boundary nodes may then vary in an open subset of their admissible boundary face. 

 Hence, locally near $\mathcal{T}\in \mathcal{D}_{\mathcal{C},L}$, $d\mu_{N,C}=d\mathcal{H}_{r}d\mathrm{m}_{\mathbb{R}^{(|\mathcal{N}(\mathcal{T})|-L)d}}$, where $\mathrm{m}$ denotes the Lebesgue measure and $\mathcal{H}_r$ denotes the $r$-Hausdorff measure, where $r$ is the Hausdorff dimension of the set of allowed boundary node locations.
  Thus, by Fubini's theorem
\begin{align*}
\textup{vol}(\mathcal{B}_{\mathcal{C},L})\leq \int_{(\partial\Omega)^L } \textup{vol}(\mathcal{B}_{\mathcal{C},L}\cap \{\mathcal{T}\,:\, \mathcal{Z}_{\mathcal{T
}}=\{x_j\}_{j=1}^L\}) d\mathcal{H}_r=0.
\end{align*}
 \end{proof}

\section*{Acknowledgements}

 JG, MM, and ES were supported by ERC Synergy Grant ``PSINumScat" 101167139.
JG was supported by EPSRC grants EP/V001760/1 and EP/V051636/1 and  the Leverhulme Trust under Research Project Grant  RPG-2023-325. 

\appendix

\section{Definition of conforming triangulations}\label{sec:App_definitions}
\begin{definition}[Conforming triangulations]
\label{d:conformingMesh}
A finite collection of nonempty $d$-simplices (see Definition~\ref{d:msimplex}) $\mathcal{T}$ is a \emph{conforming triangulation} of $\Omega\subset \mathbb{R}^d$ if
\begin{itemize}
\item[(i)] for each $T \in \mathcal{T}$, $T$ is closed, $T^\circ$ is nonempty and connected,
\item[(ii)] $\overline{\Omega}= \cup_{T\in \mathcal{T}}T$
\item[(iii)] if $T_1, T_2\in \mathcal{T}$ and $T_1^\circ\cap T_2^\circ \neq \emptyset$ then $T_1=T_2$.
\item[(iv)] for each $T \in \mathcal{T}$, $T$ has Lipschitz boundary and 
$$
\partial T=\bigcup_{j=1}^d \bigcup_{i=1}^{N_{j,d}}F_{j,i}(T)
$$
with $F_{j,i}(T)$ an open, embedded simplex of dimension $d-j$ such that $\partial \big(F_{j,i}(T)\big)=\partial\overline{F_{j,i}(T)}$,
\item[(v)] if $T_1,T_2\in\mathcal{T}$ and $F_{j_1,i_1}(T_1)\cap F_{j_2,i_2}(T_2)\neq \emptyset$, then $F_{j_1,i_1}(T_1)=F_{j_2,i_2}(T_2)$. 
\end{itemize}
\end{definition}

\begin{definition}[$m$-simplex]\label{d:msimplex}
Let $m>0$. A set $T\Subset \mathbb{R}^d$ is an \emph{$m$-simplex} if there are $\{x_j\}_{j=0}^m\subset \mathbb{R}^d$ such that
$$
T=\text{convex hull}\big(\{ x_j\}_{j=0}^m\big)
$$
$T$ is an \emph{open $m$-simplex} if $T=(\widetilde{T})^\circ$ where $\widetilde{T}$ is an $m$-simplex  and $T\neq \emptyset$. 
\end{definition}

\bibliographystyle{siamplain}
\bibliography{references}

\end{document}